\documentclass[12pt,a4paper,oneside]{article}
\usepackage{amsmath,amstext,amssymb,amsopn,amsthm,color}
\theoremstyle{plain}
\newtheorem{thm}{Theorem}[section]

\newtheorem{lem}[thm]{Lemma}

\newtheorem{cor}[thm]{Corollary}

\theoremstyle{definition}

\theoremstyle{remark}
\newtheorem*{rem*}{Remark}

\newcommand{\R}{\mathbb{R}}

\newcommand{\C}{\mathbb{C}}
\renewcommand{\H}{\mathbb{H}}

\renewcommand{\leq}{\leqslant}
\renewcommand{\geq}{\geqslant}

\renewcommand{\leq}{\leqslant}
\renewcommand{\geq}{\geqslant}

\newcommand{\pref}[1]{(\ref{#1})}

\def\o{\over}
\def\({\left(}
\def\){\right)}
\def\[{\left[}
\def\]{\right]}
\def\<{\langle}
\def\>{\rangle}

\title {Bessel Potentials, Hitting Distributions and Green Functions
\footnotetext{2000 MS Classification:
    Primary 60J65; Secondary 60J60.
    {\it Key words and phrases}: Bessel potentials, Riesz kernels,
    relativistic process, stable process, Poisson kernel, Green function,
    half-spaces.  Research partially supported by KBN
    grant 1 P03A 020 28 and RTN Harmonic Analysis and Related Problems
    contract HPRN-CT-2001-00273-HARP}}
\author{ T. Byczkowski, M. Ryznar and J. Ma{\l}ecki\\
 Institute of Mathematics and Computer Sciences,\\
  Wroc\l{}aw University of Technology, Poland}
\date{}

\begin{document}
\maketitle

\begin{abstract} The purpose of this paper is to find explicit formulas for basic objects
pertaining the local potential theory of the operator $(I-\Delta)^{\alpha/2}$,
$0<\alpha<2$. The potential theory of this operator is based on Bessel potentials
$J_{\alpha}=(I-\Delta)^{-\alpha/2}$.
We compute the {\it harmonic measure} of the half-space and write a concise form
 of the corresponding {\it Green function} for the operator $(I-\Delta)^{\alpha/2}$.
To achieve this we analyze the so-called {\it relativistic $\alpha$-stable process}
 on $\R^d$ space, killed when exiting the half-space. In terms of this process we
are dealing here with the $1$-{\it potential theory} or, equivalently, potential
theory of Schr{\"o}dinger operator based on the generator of the process with Kato's
potential $q=-1$.

\end{abstract}

\newpage
%%%%%%%%%%%%%%%%%%%%%%%%%%%%%%%%%%%%%%%%%%%%%%%%%%%%%%%%%%%%%%%%%%%%%%%%%%%%%%%%%%%%%%%%%%%%
%%%%%%%%%%%%%%%%%%%%%%%%%%%%%%%%%%%%%%%%%%%%%%%%%%%%%%%%%%%%%%%%%%%%%%%%%%%%%%%%%%%%%%%%%%%%%
\section{Introduction}
The so-called {\it interpolation spaces} play an important r{\^o}le in harmonic
analysis and partial differential equations (see, e.g. \cite{S} and \cite{Ho}).
The most classical ones are {\it Sobolev spaces} and {\it spaces of Bessel potentials}.
While (fractional) Sobolev spaces $ L_{\alpha}^{p}$ are defined in terms of
{\it Riesz potentials} $I_{\alpha}= (-\Delta)^{-\alpha/2}$, the latter ones employ
{\it Bessel potentials} defined as $J_{\alpha}= (I-\Delta)^{-\alpha/2}$. As Stein
pointed out
in his monograph \cite{S}, the both potentials exhibit the same {\it local}
 behaviour (as $|x| \to 0$) but the {\it global} one (as $|x| \to \infty$)
 of  $J_{\alpha}$ is much more regular. It is remarkable that both
 potentials can be analyzed in terms of stochastic processes (at least for $0<\alpha<2$);
 the Riesz potentials are closely related to $\alpha$-{\it stable rotation invariant
 L{\'e}vy process} while the Bessel potentials, in turn, can be investigated in terms of
 the so-called {\it relativistic process}.

 Potential theory based on Riesz kernels
 (or, equivalently: potential theory for $\alpha$-stable  rotation invariant
 L{\'e}vy process) is well developed and rich in explicit formulas, much like in the
 classical case of Brownian motion process. The homogeneity of Riesz kernels yields
 many elegant and transparent formulas for harmonic measure and Green function
 for such basic sets as balls and half-spaces in $\R^d$ (see e.g. \cite{BGR}). These
 formulas played an important r{\^o}le in setting up the so-called {\it boundary
 potential theory} of the operator $(-\Delta)^{\alpha/2}$ and the Schr{\"o}dinger
 operator based on it (see e.g. \cite{BB1} and \cite{BB2} or \cite{CS1}).

 In contrast to this situation, up to now, there were no explicit formulas known either
 for harmonic measure or  Green function for the relativistic process for sets such
 as half-planes or balls. Nevertheless, an adequate boundary potential theory (at least
 for bounded smooth sets) was set up by Ryznar \cite{Ry}. Let us point out that in recent
 years a number of publications concerning the potential theory
  of relativistic process and in particular Green function of this process
 appeared (\cite{CS2},
 \cite{K}, \cite{KL}, \cite{GRy}); all however restricted to bounded sets.

 In this paper we provide the explicit formulas for harmonic measure and Green function
 for half-spaces for the operator $-(I-\Delta)^{\alpha/2}$. Our approach consists of
 considering the operator $H_{\alpha}=I-(I-\Delta)^{\alpha/2}$, which is the infinitesimal
 generator of the relativistic $\alpha$-stable process and regarding the operator
  $-(I-\Delta)^{\alpha/2}$ as the Schr{\"o}dinger operator of the form
  ${\cal{S}}^{\alpha}=H_{\alpha}+qI$, with $q \equiv -1$. We then identify Bessel kernels
  as $1$-potentials for the relativistic $\alpha$-stable process and compute
  the corresponding $1$-Poisson kernel and $1$-Green function for half-spaces.
 It is remarkable and surprising that even though our kernels do not exhibit either
 form of homogeneity, the resulting formulas are very transparent and very similar to
 these for $\alpha$-stable case.

 The organization of the paper is as follows.
 We compute first the formulas for harmonic measure and Green function for one-dimensional
 case. To do this we apply complex-variables methods and some real-variables
 manipulations with definite integrals to obtain a satisfactory form of Green function.
 The $d$-dimensional case is then settled via application of $(d-1)$-dimensional
 Fourier transform. For technical reasons, we have to consider the Poisson kernel and
 Green function not only for the operator $-(I-\Delta)^{\alpha/2}$ but also for
 $-(m^{2/\alpha}I-\Delta)^{\alpha/2}$.
 Let us point out that we do not apply Kelvin's transform (which is
 the indispensable tool in  the $\alpha$-stable case). When performing a suitable
 limiting procedure we obtain the well-known formulas for the $\alpha$-stable case.
 The last section is devoted to various estimates for the Green function of the
 half-space $\H$, computed for the relativistic process (that is, for the operator
 $H_{\alpha}=I-(I-\Delta)^{\alpha/2}$).
 To distinguish it from the corresponding
 object computed for the operator $-(I-\Delta)^{\alpha/2}$ we call it
 {\it 0-Green function}.  The estimates are precise for $x, y \in \H$ such that
  $|x-y|<1$. To the best of our knowledge, it is the first case when the Poisson kernel
  and Green function
 of the relativistic $\alpha$-stable process of an unbounded set are examined.

%%%%%%%%%%%%%%%%%%%%%%%%%%%%%%%%%%%%%%%%%%%%%%%%%%%%%%%%%%%%%%%%%%%%%%%%%%%%%%%%%%%%%%%%%%%
%%%%%%%%%%%%%%%%%%%%%%%%%%%%%%%%%%%%%%%%%%%%%%%%%%%%%%%%%%%%%%%%%%%%%%%%%%%%%%%%%%%%%%%%

\section{Preliminaries}

    Throughout the paper  by $c, C, C_1\,\dots$ we  denote
       nonnegative constants which may depend on other constant parameters only.
        The value of $c$ or $C, C_1\,\ldots$ may change from line to line in a chain
         of estimates.

      The notion $p(u)\approx q(u),\ u \in A$ means that the ratio
      $p(u)/ q(u),\ u \in A$ is bounded from  below and above
      by positive constants which may depend on other constant parameters only.

We present in this section some basic material regarding the $\alpha$-stable
relativistic process. For more  detailed information, see
\cite{Ry} and \cite{C}. For questions regarding Markov and strong Markov property, semigroup
properties, Schr\"odinger operators and basic potential theory, the reader is referred
to \cite{ChZ} and \cite{BG}.

We first introduce an appropriate class of subordinating processes.
 Let $\theta_{\alpha}(t,u)$, $u, t>0$, denote the density function
of the strictly $\alpha/2$-stable positive standard subordinator with the Laplace
transform
$ e^{-t \lambda^{\alpha/2}}\,, \quad 0<\alpha<2\,$.

Now for $m>0$ we define another subordinating process $T_{\alpha}(t,m)$
modifying the corresponding probability density function in the following way:

 \begin{equation*}
 \theta_{\alpha}(t,u,m)= e^{mt}\,\theta_{\alpha}(t,u)\,e^{-m^{2/\alpha}u},  \quad u>0\,.
 \end{equation*}
 We derive the Laplace transform of $T_{\alpha}(t,m)$ as follows:
 \begin{equation} \label{subord2}
 E^{0} e^{-\lambda T_{\alpha}(t,m)} =e^{mt}\, e^{-t (\lambda+m^{2/\alpha})^{\alpha/2}}\,.
 \end{equation}
 Let $B_t$ be the symmetric Brownian motion in $\R^d$ with the characteristic
  function of the form
 \begin{equation} \label{brownian}
 E^{0}e^{i\xi \cdot B_t} = e^{-t|\xi|^2}\,.
 \end{equation}
 Assume that the processes $T_{\alpha}(t,m)$ and $B(t)$ are stochastically independent.
 Then the process $X_t^{\alpha,m}= B_{T_{\alpha}(t,m)}$ is called the
  $\alpha$-stable relativistic process
 (with parameter $m$).
  In the sequel we use the generic notation $X_t^m$ instead of
   $X_t^{\alpha,m}$.
  If $m=1$ we write
  $T_{\alpha}(t)$ instead of $T_{\alpha}(t,m)$ and $X_t$ instead of $X_t^1$.

  $X_t^m$ is a L\'evy process (i.e. homogeneous, with independent
 increments). We always assume that sample paths of the process $X_t^m$ are right-continuous
 and have left-hand limits ("cadlag"). Then $X_t^m$ is Markov and has the strong
 Markov property under the so-called standard filtration.

 Various potential-theoretic objects in the theory of the process $X_t^m$ are expressed
 in terms of modified Bessel functions $K_{\nu}$ of the second kind, called also Macdonald
 functions. For convenience of the reader
 we collect here basic informations about these functions.
%%%%%%%%%%%%%%%%%%%%%%%%%%%%%%%%%%%%%%

  $K_\nu,\ \nu\in \R $, the modified Bessel function of the second kind
   with index $\nu$,  is given by the
 following formula:
 $$ K_\nu(r)=  2^{-1-\nu}r^\nu \int_0^\infty e^{-v}e^{- {r^2\o 4v}}v^{-1-\nu}dv\,,
 \quad r>0.$$
  For properties of  $K_\nu$ we refer the reader to \cite{E1}. In the sequel we will
  use the asymptotic behaviour of $K_\nu$:

 \begin{eqnarray}
  K_\nu(r)& \cong& {\Gamma(\nu)\o 2} \left({r\o 2}\right)^{\nu}\,\quad  r\to 0^+, \quad \nu>0,
  \label{asympt0}\\
    K_0(r)&\cong& -\log r,  \quad r\to 0^+, \label{asympt00}\\
    K_\nu(r)&\cong& \frac {\sqrt{\pi}} {\sqrt{2 r}}\,e^{-r},  \quad r\to \infty, \label{asympt_infty}
 \end{eqnarray}
 where  $g(r) \cong f(r) $ denotes that the ratio of $g$ and $f$ tends to $1$.
  For  $\nu <0 $ we have $K_\nu(r)=K_{-\nu}(r)$, which determines the asymptotic
  behaviour for negative indices.

%%%%%%%%%%%%%%%%%%%%%%%%%%%%%%%%%%%%%%%%

  The  $\alpha$-stable relativistic  density (with parameter $m$)
   can now be computed in the following way:
  \begin{equation} \label{reldensity0}
 p^m_t(x)=\int_0^\infty\theta_{\alpha}(t,u,m)\, g_u(x) du,
 \end{equation}
 where $g_u(x)$ is  the Brownian semigroup, defined by \pref{brownian}. A particular case
  when $\alpha=1$ is called
 the relativistic Cauchy semigroup on $\R^{d}$ with parameter $m$.
 The formula below exhibits the explicit form of this density:
  \begin{lem}[relativistic Cauchy semigroup]
   The density $\tilde{p}^m_{t}$ of the relativistic Cauchy process is of the form:
 \begin{equation} \label{Cauchyrel}
   \tilde{p}^m_{t}(x)=
    2(m/2\pi)^{(d+1)/2}\, te^{mt}
     {K_{(d+1)/2}(m(|x|^2+t^2)^{1/2})\o (|x|^2+t^2)^{d+1 \o 4}}.
   \end{equation}
  \end{lem}
  \begin{proof}
  We carry out the corresponding computations as follows:
  \begin{eqnarray*}
    \tilde{p}^m_{t}(x)&=&
    e^{mt} \int_0^{\infty} {1 \o (4\pi u)^{d/2} } e^{-|x|^2/4u} e^{-m^2u}
    {t \o \sqrt{4\pi}} u^{-3/2} e^{-t^2/4u}\, du \\
    &=& {t e^{mt} \o (4\pi)^{d+1 \o 2}} \int_0^{\infty} e^{-m^2 u}
    e^{-(|x|^2+t^2)/4u} {du \o u^{{d+1 \o 2}+1}} \\
    &=&    2(m/2\pi)^{(d+1)/2}\, te^{mt}
  {K_{(d+1)/2}(m(|x|^2+t^2)^{1/2})\o (|x|^2+t^2)^{d+1 \o 4}}.
  \end{eqnarray*}
 \end{proof}

In the case of arbitrary $\alpha$, $0<\alpha <2$ we have the following useful estimate
(see \cite{Ry}):

  \begin{lem} \label{transden_0}
  There exists a constant $c=c(\alpha ,d,m)$ such that
  \begin{equation} \label{transupper}
  \max_{x \in {\R}^d} p_t^m(x) \le c (t^{-d/2} + t^{-d/{\alpha}}) \,.
  \end{equation}
  \end{lem}

 In the next lemma we compute the  Fourier transform of the transition
 density \pref{reldensity0}:
 \begin{lem}[Fourier transform of $p^m_t$]
 The Fourier transform of $\alpha$-stable relativistic density $p^m_t$is of the form:
 \begin{equation} \label{Fourierrel}
\widehat{p^m_t}(z) =
  e^{mt}e^{-t(|z|^2+m^{2/\alpha})^{\alpha/2}}.
 \end{equation}
 \end{lem}
 \begin{proof}
 \begin{eqnarray*}
 \widehat{p^m_t}(z)&=&\int_{\R^{d}}  p^m_t(x)e^{i(z,x)}dx
 = \int_{\R^{d}}\int_0^{\infty} e^{mt} g_u(x)e^{-m^{2/\alpha}u}\theta_{\alpha}(t,u)du\,
  e^{i(z,x)}dx \\
   &=& e^{mt} \int_0^{\infty}e^{-u|z|^{2}}e^{-m^{2/\alpha}u}\theta_{\alpha}(t,u)du
   = e^{mt} \int_0^{\infty}e^{-u(|z|^{2}+m^{2/\alpha})}\theta_{\alpha}(t,u)du \\
   & = & e^{mt}e^{-t(|z|^2+m^{2/\alpha})^{\alpha/2}}.
   \end{eqnarray*}
 \end{proof}
  Specifying this to the case $\alpha=1$ we obtain
 \begin{equation}\label{relFourier}
  \widehat{\tilde{p}^m_t}(z)=
 e^{mt}e^{-t(|z|^2+m^{2})^{1/2}}.
    \end{equation}

 From the form of the Fourier transform we have the following scaling property:
 \begin{equation}\label{scaling}
 {p}^m_{t}(x)=m^{d/\alpha}{p}^1_{mt}(m^{1/\alpha}x).
   \end{equation}
 In terms of one-dimensional distributions of the relativistic process (starting from
 the point $0$) we obtain
 \begin{equation*}
 X_t^m \sim m^{-1/\alpha} X_{mt}^1\,,
 \end{equation*}
 where $X_t^m$ denotes the relativistic $\alpha$-stable process with parameter $m$
 and "$\sim$" denotes equality of distributions.
 Because of this scaling property, we usually restrict our attention to the case
  when $m=1$, if not specified otherwise.
  When $m=1$ we omit the superscript "$1$", i.e. we write $p_{t}(x)$
  instead of $p_t^1(x)$.

 We work here within the framework of the so-called $\lambda$-{\it potential theory},
 for $\lambda>0$. The standard reference book on general potential theory is the
 monograph \cite{BG}. For
 convenience of the reader we collect here the basic information with emphasis on
 what is known (and needed further on) about the $\alpha$-stable relativistic process.

 We begin with the kernel of the resolvent semigroup of the process $X_t^m$. We call
 the kernel corresponding to parameter $\lambda$ the $\lambda$-{\it potential}
 of the process and denote it by $U^m_{\lambda}(x)$ . It is of a particular
 simple form when $\lambda=m$ and we state it for further references.
  Again we
  denote by $U_{\lambda}(x)$ the $\lambda$-potential for the case $m=1$.

  \begin{lem} [$m$-potential for relativistic process with parameter $m$]
  \begin{equation}\label{m-potential}
  U^m_m(x)=C(\alpha,d)\,m^{d-\alpha \o 2 \alpha}\,
   {K_{(d-\alpha)/2}(m^{1/\alpha}|x|)\o  |x|^{(d-\alpha)/2}}\,,
  \end{equation}
  where $C(\alpha,d)= {2^{1-(d+\alpha)/2} \o {\Gamma(\alpha/2)\pi^{d/2}}}$.
  \end{lem}
  \begin{proof}
  We provide  calculations for $m=1$ and the general case follows from (\ref{scaling}).
  \begin{eqnarray*}
    U_1(x)&=&\int_0^\infty e^{-t} p_t(x)dt=\int_0^\infty
   \int_0^\infty g_u(x)e^{-u}\theta_{\alpha}(t,u)du\,dt \nonumber\\
    &=& \int_0^\infty {1 \o {(4\pi u)^{d/2}}}
    e^{-\frac {|x|^2}{4u}}e^{-u} \left(\int_0^\infty\theta_{\alpha}(t,u)dt\right)du \nonumber\\
    &=&{1 \o {\Gamma(\alpha/2)(4\pi)^{d/2}}}\int_0^\infty
     e^{-\frac {|x|^2}{4u}}e^{-u} {du \o u^{{d-\alpha \o 2} +1}} \nonumber\\ &=&
    {2^{1-(d+\alpha)/2} \o {\Gamma(\alpha/2)\pi^{d/2}}}
    { K_{(d-\alpha)/2})\o |x|^{(d-\alpha)/2}}  \,.
    \end{eqnarray*}

    \end{proof}

    The {\it first exit time} of an (open)
   set  $D\subset {\R}^{d}$
   by the process $X_t^m$ is defined by the formula
   \begin{equation*}
   \tau_{D}=\inf\{t\geq 0;\, X_t^m\notin D\}\,.
   \end{equation*}

  The basic object in potential theory of $X_t^m$ is the
  $\lambda$-{\it harmonic measure}  of the
  set $D$. It is defined by the formula:

  \begin{equation} \label{harm_def}
  P_D^{\lambda}(x,A)=
  E^x[\tau_D<\infty; e^{-\lambda \tau_D} {\bf{1}}_A(X_{\tau_D}^m)].
  \end{equation}
  The density kernel of  the measure $P_D^{\lambda}(x,A)$ (if it exists) is called the
  $\lambda$-{\it Poisson kernel} of the set $D$.

  In Section 3 and 4 we prove the existence of the $m$-Poisson kernel providing at
    the same time an explicit formula for it. The more general case of
     existence of $\lambda$-Poisson kernel
    can be deduced from papers \cite{IW} and \cite{BS}.

  Another fundamental object of potential theory is the {\it killed process} $X_t^{m,D}$
  when exiting the set $D$. It is defined in terms of sample paths up to time $\tau_D$.
  More precisely, we have the following "change of variables" formula:
  \begin{equation} \label{transkilled_den}
  E^x f(X_t^{m,D}) =  E^x[t<\tau_D; f(X_t^{m})]\,, t>0\,.
  \end{equation}
  The density function of transition probability of the process $X_t^{m,D}$ is denoted
  by $p_t^{m,D}$. We have
  \begin{equation*}
  p_t^{m,D}(x,y) = p_t^{m}(x-y) -
   E^x[t\ge \tau_D; p_{t-\tau_D}^m(X_{\tau_D}^m-y)]    \,, \quad x, y \in {\R}^d\,.
  \end{equation*}
  Obviously, we obtain
   \begin{equation*}
     p_t^{m,D}(x,y) \le p_t^{m}(x,y) \,, \quad x, y \in {\R}^d\,.
   \end{equation*}

  $ p_t^{m,D}$ is a strongly contractive semigroup (under composition) and shares most
  of properties of the semigroup $ p_t^{m}$. In particular, it is strongly Feller and
  symmetric: $ p_t^{m,D}(x,y) =  p_t^{m,D}(y,x)$.  When $m=1$, we write, as before,
   $ p_t^{D}$, instead of $ p_t^{1,D}$.

  The $\lambda$-potential of the process $X_t^{m,D}$ is called the
  $\lambda$-{\it Green function} and is denoted by $G_D^{\lambda}$. Thus, we have
  \begin{equation*}
   G_D^{\lambda}(x,y)= \int_0^{\infty} e^{-\lambda t}\,p_t^{m,D}(x,y)\,dt\,.
  \end{equation*}

  The "sweeping out" procedure provides another formula for the $\lambda$-{\it Green
  function} of the set $D$ in terms of the $\lambda$-harmonic measure:
  \begin{equation} \label{green_def}
  G_D^{\lambda}(x,y)= U_{\lambda}^m -
  \int_{{\R}^d} U_{\lambda}^m(z,y)\, P_D^{\lambda}(x,dz) \,.
  \end{equation}

  We now state some basic scaling properties both for the $m$-Poisson kernel and
  the $m$-Green function. The proof employs the scaling property \pref{scaling} and
  consists of elementary but tedious calculation and is omitted.

  \begin{lem}[Scaling Property] \label{scaling1}
Let $D$ be an open subset of $\R^d$ and $P_D^{m}$,  $G_D^{m}$ be $m$-Poisson kernel,
or $m$-Green function, respectively, for $D$. Then

$$P_D^{m}(x,u)= m^{d/\alpha}P^1_{m^{1/\alpha}D}(m^{1/\alpha}x,m^{1/\alpha}u),\quad
x\in D, u\in D^c\,, $$

$$G_D^{m}(x,y)= m^{(d-\alpha)/\alpha}G^1_{m^{1/\alpha}D}(m^{1/\alpha}x,m^{1/\alpha}y),\quad x\in D, y\in D\,. $$
Thus, if $D$ is a cone with vertex at $0$ we obtain:
$$P_D^{m}(x,u)= m^{d/\alpha}P^1_{D}(m^{1/\alpha}x,m^{1/\alpha}u),\quad x\in D, u\in D^c\,, $$

$$G_D^{m}(x,y)= m^{(d-\alpha)/\alpha}G^1_{D}(m^{1/\alpha}x,m^{1/\alpha}y), \quad x\in D,
 y\in D\,. $$
   \end{lem}

 Taking into account the above lemma and the fact that we restrict our attention
 to the case when $D=\H$, the upper half-space, we usually assume that
 $\lambda=m=1$ and, as before, omit superscript $"1"$, that is, we write
 $P_D(x,u)$ and $G_D(x,y)$ instead of
  $P_D^{1}(x,u)$ and $G_D^{1}(x,y)$.

 We also recall the form of the density function $\nu^m(x)$ of the L{\'e}vy measure of the
 relativistic $\alpha$-stable process (see \cite{Ry}):
  \begin{lem}[L{\'e}vy measure of relativistic process with parameter $m$]
  \begin{equation} \label{levymeasure}
  \nu^m(x)= {\alpha 2^{\alpha -d \o 2} \o \pi^{d/2} \Gamma(1- {\alpha \o 2})}
  \left( {m^{1/\alpha} \o |x|} \right)^{ d+\alpha \o 2}
  K_{ d+\alpha \o 2}(m^{1/\alpha}|x|)\,.
  \end{equation}
  \end{lem}

 When $m=1$ we omit, as usual, the superscript "$1$".

    The next lemma is taken from  \cite{G}. For reader's convenience we
     provide its proof.
     \begin{lem}    \label{Grzywny}
     For $m=1$ we have
\begin{equation} \label{Grzywny1}P^x(\tau_\H\ge t)\le
 C \frac{x_d + \ln t}{t^{1/2}}\,,\quad t\ge 2\,,\, x_d>0.\end{equation}
\end{lem}

\begin{proof}
Let $Y_t=X^{(d)}_t$, where $X_t=(X^{(1)}_t,\ldots,X^{(d)}_t)$. By the symmetry of the
random variable $Y_t$ we obtain
$$P^x (\tau_\H>t)= P^x (\inf_{s\leq t}Y_s >0)=  P^0 (\inf_{s\leq t}(-Y_s+x_d) >0)=
P^0 (\sup_{s\leq t}Y_s < x_d).$$
Using a version of the L\'evy inequality (\cite{B}, Ch.7, 37.9) we have for any
 $\varepsilon,y>0$ that
$$
2 P^0 (Y_t\geq y+2\varepsilon)-2\sum^n_{k=0}P^0(Y_{\frac{tk}{n}}-Y_{\frac{t(k-1)}{n}}\geq \varepsilon)
\leq P^0 (\sup_{k\leq n}Y_{\frac{tk}{n}} \geq y).
$$
Note that  $\sum^n_{k=0}P^0(Y_{\frac{tk}{n}}-Y_{\frac{t(k-1)}{n}}\geq \varepsilon)=nP^0(Y_{\frac{t}{n}}\geq \varepsilon)\to t\int^{\infty}_\varepsilon\nu(x)dx$,
hence, by symmetry again
$$ P^0 (|Y_t|\geq y+2\varepsilon)-2t\int^{\infty}_{\varepsilon}\nu(x)dx =
2 P^0 (Y_t\geq y+2\varepsilon)-2t\int^{\infty}_{\varepsilon}\nu(x)dx\leq
P^0 (\sup_{s\leq t}Y_t \geq y).$$
This implies that
$$P^x (\tau_\H>t)=P^0 (\sup_{s\leq t}Y_s < x_d)\le P^0 (|Y_t| <x_d+2\varepsilon)
+2t\int^{\infty}_{\varepsilon}\nu(x)dx\,. $$
For $\varepsilon\geq 1$ we obtain from \pref{levymeasure} and \pref{asympt_infty}
$$\int^{\infty}_{\varepsilon}\nu(x)dx\leq C e^{-\varepsilon}\varepsilon^{-\alpha/2}.$$
Because of the Lemma \ref{transden_0}, we have
 that the density of $Y(t)$ is bounded by $Ct^{-1/2}, t\ge 2$ hence taking
$\varepsilon=\frac{3}{2}\ln t$ we obtain
$$P^x (\tau_\H>t)\leq
 C \left(x_d+\ln t\right) t^{-1/2}.$$

\end{proof}

    In order to improve the above estimate for  $x$ close to the boundary we use the
    following result proved recently in \cite{GRy}.
    \begin{lem}    \label{GrzywnyRyznar} Assume that d=1. Let $D=(0,1)$ and  $x\in D$. Then
     $$P^x ( X_{\tau_D}\in dz)\approx  \frac {(x(1-x))^{\alpha/2}}{(z-1)^{\alpha/2}(z-x)}\ e^{-z}
     , \quad x\in D,\ z>1.$$
     This implies that
    $$ E^x [X_{\tau_D}>1;X_{\tau_D}]\approx P^x ( X_{\tau_D}>1)\approx x^{\alpha/2}.$$
    %%%%%%%%%%%%%%
    We also have that $$E^x\tau_D \approx {(x(1-x))^{\alpha/2}}.$$
    \end{lem}
    Actually in \cite{GRy} it was shown that the Green function of $D$ is comparable with
     the Green function of the corresponding stable process. By standard arguments
      (see \cite{Ry}) this implies the above lemma.

    %%%%%%%%%%%%%%%%
    We further need the following strenghtening of Lemma \ref{Grzywny}.

    \begin{lem}  \label{Grzywnyimprove}
     For $0<x_d<1/2$ we have
     \begin{equation}\label{uboundtail}
     P^x(\tau_\H\ge t)\le C x_d^{\alpha/2}\;\ln t/t^{1/2},\quad t\ge 2,
     \end{equation}
     where $C$ is a constant.
     \end{lem}
     \begin{proof}
     It is enough to prove the claim for $d=1$. Let $D=(0,1)$ and assume that $0<x<1/2$.

     By the Markov property and then by  Lemma \ref{Grzywny} we obtain for $t\ge 2$:

     \begin{eqnarray}
     P^x(\tau_\H\ge t)&=& P^x(\tau_D\ge t, \tau_D=\tau_\H)+
     E^x[\tau_D<\tau_\H\,; P^{X_{\tau_D}}(\tau_\H\ge t]\nonumber\\
     &\le& P^x(\tau_D\ge t)+
     E^x [\tau_D<\tau_\H; X_{\tau_D} + \ln t]/t^{1/2}\nonumber\\
     &=& P^x(\tau_D\ge t)+
     E^x [X_{\tau_D}>1\,;X_{\tau_D}]/t^{1/2}+\ln t\;P^x (X_{\tau_D}>1)/t^{1/2}\nonumber\\
     &\le&C x^{\alpha/2}\ln t/t^{1/2}.\nonumber
    \end{eqnarray}
    %%%%%%%%%%%%%%%%%%%%%%%%%
     The last inequality follows from Lemma \ref{GrzywnyRyznar}. The proof is complete.
     \end{proof}

%%%%%%%%%%%%%%%%%%%%%%%%%%%%%%%%%%%%%%%%%%%%%%%%%%%%%%%%%%%%%%%%%%%%%%%%%%%%%%%%%%%%%%%%%%%
%%%%%%%%%%%%%%%%%%%%%%%%%%%%%%%%%%%%%%%%%%%%%%%%%%%%%%%%%%%%%%%%%%%%%%%%%%%%%%%%%%%%%%%%%

  \section{Poisson kernel and Green function of
  $-(m^{2/\alpha}I-{d^2 \o dx^2})^{\alpha/2}$ for
   $\R^{+}$ }
  We first consider the one-dimensional case. For technical reasons (see Section 4)
  we have to distinguish notation in this case from the $d$-dimensional one so
  we denote the relativistic Poisson kernel by  $ Q_{(0,\infty)}^m $ (instead of
  $P_{\H}^m $ as previously).
  \begin{thm}[one-dimensional relativistic  Poisson kernel]\label{1dPoiss}
  Denote
  \begin{equation*}
  E^x[e^{-m\tau_{(0,\infty)}}; X^m_{\tau_{(0,\infty)}} \in du] =
   Q_{(0,\infty)}^m(x,u)\,.
  \end{equation*}
  Then we have for $u<0<x$:

  \begin{equation*}
 Q_{(0,\infty)}^m(x,u) =
 \frac{\sin(\pi\alpha /2)}{\pi}
 \left(\frac{x}{-u}\right)^{\alpha /2}\frac{e^{-m^{1/\alpha}(x-u)}}{x-u}\,.
 \end{equation*}
 \end{thm}
  \begin{proof} By the scaling property it is enough to deal with $m=1$. According to
  the general facts of potential theory (see e.g. \cite{BGR})
  it is enough to verify that the measure $Q_{(0,\infty)}(x,u)$ satisfies the following
  "sweeping out" principle:
  \begin{equation*}
  \int_{-\infty}^0 Q_{(0,\infty)}(x,u) U_1(u,y)\,du =  U_1(x,y)\,,
  \end{equation*}
  for all $y<0$, where $U_1(x,y)$ is the $1$-potential of the relativistic
   $\alpha$-stable process (see \pref{m-potential}). This is, however,
    shown in the following lemma below.
  \end{proof}

  \begin{lem}
  For $x>0>y$ we have
  \begin{eqnarray*}
  \frac{\sin({\pi\alpha\o2})}{\pi}\int_{-\infty}^0
  \left(\frac{x}{-u}\right)^{\alpha\o2}\frac{e^{-|x-u|}}{|x-u|}\frac{K_{1-\alpha\o2}(|u-y|)}{|u-y|^{1-\alpha\o2}}\,du
  =\frac{K_{1-\alpha\o2}(|x-y|)}{|x-y|^{1-\alpha\o2}}
  \end{eqnarray*}
  \end{lem}
  \begin{proof}
  Let $x>0>y$ and consider the following function of complex
  variable $z$:
  \begin{equation*}
  f(z) =
  \frac{1}{z^{\alpha\o2}}\frac{e^{z-x}}{x-z}\frac{K_{1-\alpha\o2}(z-y)}{(z-y)^{1-\alpha\o2}}.
  \end{equation*}
  This function is holomorphic in $\C\backslash(-\infty,0]
  \backslash\{x\}$. We are going to integrate the function above over
  a contour described below. We make the branch cut along the axis
  $(-\infty,0]$ and make the contour of integration to wrap
  around this line (see picture) and we add a circle around $x$ say
  $\gamma_3$.

  
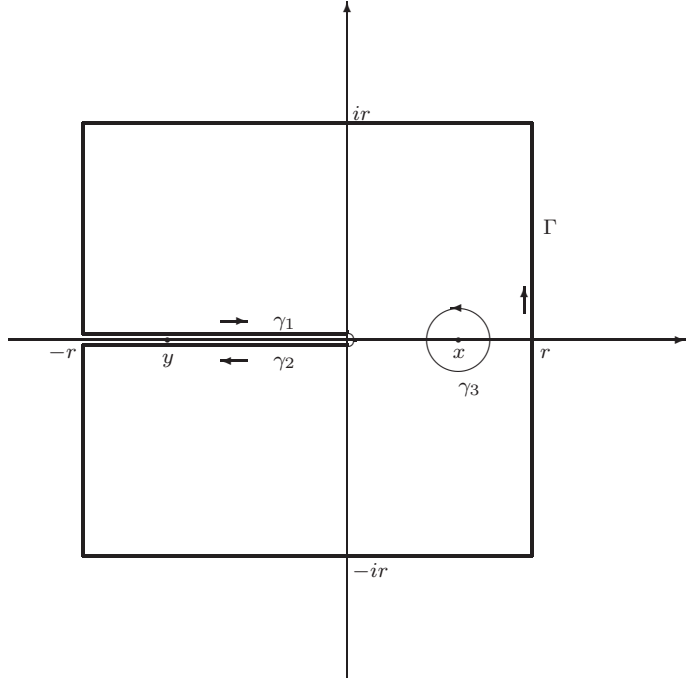
\begin{center}
    \begin{picture}(256,256)(0,0)
      \put(0,128){\vector(1,0){256}}
      \put(128,0){\vector(0,1){256}}
      \linethickness{1pt}
      \put(128,128){\oval(5,5)[r]}
      \put(128,126){\line(-1,0){100}}
      \put(128,130){\line(-1,0){100}}
      \put(28,126){\line(0,-1){80}}
      \put(28,130){\line(0,1){80}}
      \put(28,210){\line(1,0){170}}
      \put(28,46){\line(1,0){170}}
      \put(198,46){\line(0,1){164}}
      \put(201,123){$_r$}\put(130,40){$_{-ir}$}
      \put(130,213){$_{ir}$}\put(15,123){$_{-r}$}
      \linethickness{0.5pt}
      \put(195,138){$\vector(0,1){10}$}
      \put(100,135){$_{\gamma_1}$} \put(100,120){$_{\gamma_2}$}
      \put(170,110){$_{\gamma_3}$}
      \put(202,170){$_{\Gamma}$}
      \put(80,135){$\vector(1,0){10}$} \put(90,120){$\vector(-1,0){10}$}
      \put(168,140){$\vector(-1,0){1}$}
      \put(60,128){\circle*{2}}\put(170,128){\circle*{2}}
      \put(168,123){$_x$}\put(58,121){$_y$}
      \put(170,128){\circle{25}}
     \end{picture}
     \end{center}
   By the Cauchy theorem, we get
  \begin{equation}
  \label{cauchytheorem} \frac{1}{2\pi i}\int_\Gamma
  f(z)dz+\frac{1}{2\pi
  i}\left(\int_{\gamma_1}+\int_{\gamma_2}\right)f(z)dz=\frac{1}{2\pi
  i} \int_{\gamma_3}f(z)dz
  \end{equation}
  Using the asymptotic expansion for the modified Bessel function
  $K_v(z)$ we get
  \begin{equation*}
      \frac{K_{1-\alpha\o2}(z-y)}{(z-y)^{1-\alpha\o2}} =
      \frac{e^{y-z}}{(z-y)^{1-{\alpha\o2}}}R_0(z-y)\/,
  \end{equation*}
  where $R_0(z)=O(1)$ and $z$ is large enough. Using the expression
  given above it is easy to show that the integral over $\Gamma$ of
  $f(z)$ vanishes when $r \to \infty$.

  To calculate the integrals over $\gamma_1$ and $\gamma_2$ we examine
  the behaviour of the function $f$ near the branch cut. We have the
  following relations for the modified Bessel functions for $u<y$
  \begin{eqnarray*}
      \lim_{\epsilon \to 0+}K_{1-\alpha\o2}(u-y+i\epsilon) &=&
      e^{{i\pi(\alpha-1)\o2}}K_{1-\alpha\o2}(|u-y|)-i\pi I_{1-\alpha\o2}(|u-y|)\/,\\
      \lim_{\epsilon \to 0+}K_{1-\alpha\o2}(u-y-i\epsilon) &=&
      e^{{i\pi(1-\alpha)\o2}}K_{1-\alpha\o2}(|u-y|)+i\pi I_{1-\alpha\o2}(|u-y|)\/.\\
  \end{eqnarray*}
  We have also
  \begin{eqnarray*}
  \lim_{\epsilon \to 0+}(u-y+i\epsilon)^{-{1-\alpha\o2}} &=&
  e^{{i\pi(\alpha-1)\o2}} |u-y|^{-{1-\alpha\o2}}\/,\\
  \lim_{\epsilon \to 0+}(u-y-i\epsilon)^{-{1-\alpha\o2}} &=&
  e^{{i\pi(1-\alpha)\o2}} |u-y|^{-{1-\alpha\o2}}\/.
  \end{eqnarray*}
  Next, for $u<0$ we have
  \begin{eqnarray*}
  \lim_{\epsilon \to 0+}(u+i\epsilon)^{-{\alpha\o2}} &=&
  e^{-{i\pi\alpha\o2}} |u|^{-{\alpha\o2}}\/,\\
  \lim_{\epsilon \to 0+}(u-i\epsilon)^{-{\alpha\o2}} &=&
  e^{{i\pi\alpha\o2}} |u|^{-{\alpha\o2}}\/.
  \end{eqnarray*}
  Using all the relations given above we obtain that for $u<y$
  \begin{eqnarray*}
  &f_1^{+}(u) = \lim_{\epsilon \to 0+}f(u+i\epsilon)
  = -\frac{1}{|u|^{\alpha\o2}}\frac{e^{u-x}}{x-u} \left[e^{i\pi\alpha\o2}\frac{K_{1-\alpha\o2}(|u-y|)}{|u-y|^{1-\alpha\o2}}+\pi\frac{I_{1-\alpha\o2}(|u-y|)}{|u-y|^{1-\alpha\o2}}\right]\/,\\
  &f_1^{-}(u) = \lim_{\epsilon \to 0+}f(u-i\epsilon)
  = -\frac{1}{|u|^{\alpha\o2}}\frac{e^{u-x}}{x-u} \left[e^{-{i\pi\alpha\o2}}\frac{K_{1-\alpha\o2}(|u-y|)}{|u-y|^{1-\alpha\o2}}+\pi\frac{I_{1-\alpha\o2}(|u-y|)}{|u-y|^{1-\alpha\o2}}\right]\/.
  \end{eqnarray*}
  Similarly for $y<u<0$ we get
  \begin{eqnarray*}
  &f_2^{+}(u) = \lim_{\epsilon \to 0+}f(u+i\epsilon) = \frac{1}{|u|^{\alpha\o2}}\frac{e^{u-x}}{x-u}
  e^{-{i\pi\alpha\o2}}\frac{K_{1-\alpha\o2}(|u-y|)}{|u-y|^{1-\alpha\o2}}\/,\\
  &f_2^{-}(u) = \lim_{\epsilon \to 0+}f(u-i\epsilon) = \frac{1}{|u|^{\alpha\o2}}\frac{e^{u-x}}{x-u}
  e^{{i\pi\alpha\o2}}\frac{K_{1-\alpha\o2}(|u-y|)}{|u-y|^{1-\alpha\o2}}\/.
  \end{eqnarray*}
  Combining all above we find that
  \begin{eqnarray*}
  \frac{1}{2\pi i}\left(\int_{\gamma_1}+\int_{\gamma_2}\right)f(z)dz
  &\to& \frac{1}{2\pi i}
  \int_{-\infty}^{y}(f_1^{+}(u)-f_1^{-}(u))du +
  \frac{1}{2\pi i}\int_{y}^{0}(f_2^{+}(u)-f_2^{-}(u))du\\
  &=&-\frac{\sin({\pi\alpha\o2})}{\pi}\int_{-\infty}^0
  \frac{1}{(-u)^{\alpha\o2}}\frac{e^{-|x-u|}}{|x-u|}\frac{K_{1-\alpha\o2}(|u-y|)}{|u-y|^{1-\alpha\o2}}\,du\/.\\
  \end{eqnarray*}
  We have also
  \begin{eqnarray*}
  \frac{1}{2\pi i} \int_{\gamma_3}f(z)dz &=&
  -\frac{1}{x^{\alpha\o2}}\frac{K_{1-\alpha\o2}(|x-y|)}{|x-y|^{1-\alpha\o2}}\/.
  \end{eqnarray*}
  Using (\ref{cauchytheorem}) and the relations given above we obtain
  the desired formula.
  \end{proof}

  \begin{thm}[one-dimensional relativistic $\alpha$-stable Green function] \label{green_1_form}
    For $x,y>0$ we have

     \begin{eqnarray*}
  G_{(0,\infty)}^{m}(x,y) &=&
        {|x-y|^{{\alpha-1}} \o  2^\alpha \Gamma(\alpha/2)^2}
        \int_0^{{4xy\o(x-y)^2}}e^{-m^{1/\alpha}|x-y|(t+1)^{1/2}}t^{\alpha/2-1}(t+1)^{-1/2}\,dt\/.
  \end{eqnarray*}

   \end{thm}
   \begin{proof} By the  scaling property we may assume $m=1$.
    Due to the symmetry of the Green function, it is enough to determine
     $G_{(0,\infty)}(x,y)$ for $0<x<y$. We
    compute the compensator of the Green function for the one-dimensional $\alpha$-stable
     relativistic process. We want
    to evaluate the following expression
     \begin{eqnarray}
     \nonumber
         H(x,y) &=& E^x[ e^{-\tau_{(0,\infty)}} U_1(X_{\tau_{(0,\infty)}},y)] \\
     \label{Green01}
         &=&
         C(\alpha,1){\sin(\pi\alpha/2) \o  \pi}  \int_{-\infty}^0
           \left({x \o -u} \right)^{\alpha/2} {e^{-(x-u)}\o x-u}\,
             {K_{1-\alpha\o2}(y-u)\o(y-u)^{1-\alpha\o2}}\,du \/,
       \end{eqnarray}
       where
       $C(\alpha,1)={2^{(1-\alpha)/2}\o\Gamma(\alpha/2)\pi^{1/2}}$.
       Substituting $(-u)^{1-\alpha/2}=v$ and taking into account the following well-known
       identities
         \begin{eqnarray*}
          {1 \o v^{2/(2-\alpha)} + x}
           &=&  \int_0^{\infty} e^{-w(x+v^{2/(2-\alpha)})}\,dw\, , \\
           C(\alpha,1){K_{1-\alpha\o2}(z)\o z^{1-\alpha\o2}}&=& \frac{1}{\Gamma(\alpha/2)\Gamma(1-\alpha/2)}\int_1^\infty
           {e^{-zt}\o(t^2-1)^{\alpha/2}}\,dt\/,
          \end{eqnarray*}
          we obtain that the right-hand side of (\ref{Green01}) is of the form
       \begin{eqnarray*}
               \lefteqn{{2\sin(\pi\alpha/2)x^{\alpha/2}e^{-x}\o(2-\alpha)\pi\Gamma(\alpha/2)\Gamma(1-\alpha/2)}
               \int_0^\infty \Bigl\{\int_0^\infty e^{-wx}e^{-wv^{2/(2-\alpha)}}\,dw\Bigr\}
               e^{-v^{2/(2-\alpha)}}\int_1^\infty
               {e^{-(y+v^{2/(2-\alpha)})t}\o(t^2-1)^{\alpha/2}} \, dt\,dv}\\
                &=&{2\sin(\pi\alpha/2)x^{\alpha/2}e^{-x}\o(2-\alpha)\pi\Gamma(\alpha/2)\Gamma(1-\alpha/2)}
               \int_1^\infty\int_0^\infty
               e^{-wx}\int_0^\infty e^{-(w+t+1)v^{2/(2-\alpha)}}\,dv\,
               {e^{-yt}\o(t^2-1)^{\alpha/2}}\,dw\,dt\/.
               \end{eqnarray*}
          The interior integral can be expressed as
          \begin{eqnarray*}
               \int_0^\infty e^{-(w+t+1)v^{2/(2-\alpha)}}\,dv &=& (w+t+1)^{(\alpha-2)/2}\int_0^\infty
                e^{-u^{2/(2-\alpha)}}\,du\\
                &=& {2-\alpha\o2}\Gamma(1-\alpha/2)(w+t+1)^{(\alpha-2)/2}\/.
          \end{eqnarray*}
          Consequently, we get
        \begin{eqnarray*}
                 H(x,y)&=&
                 {\sin(\pi\alpha/2)x^{\alpha/2}e^{-x}\o\pi\Gamma(\alpha/2)}
                \int_1^\infty\int_0^\infty
                e^{-wx}(w+t+1)^{(\alpha-2)/2}\,dw
                {e^{-yt}\o(t^2-1)^{\alpha/2}}\,dt \\
                 &=& {\sin(\pi\alpha/2)\o\pi\Gamma(\alpha/2)}
                \int_1^\infty\left(x^{\alpha/2}\int_0^\infty
                e^{-x(w+t+1)}(w+t+1)^{(\alpha-2)/2}\,dw\right)
                {e^{-(y-x)t}\o(t^2-1)^{\alpha/2}} \, dt\/.
                 \end{eqnarray*}
                 Observe that the expression in brackets can be computed as follows
                 \begin{eqnarray*}
                 x^{\alpha/2}\int_0^\infty
                e^{-x(w+t+1)}(w+t+1)^{(\alpha-2)/2}\,dw
                 &=& \int_{x(t+1)}^\infty e^{-s}s^{{\alpha/2}-1}ds \\
                 &=&\Gamma(\alpha/2)-\int_0^{x(t+1)}
                 e^{-s}s^{{\alpha/2}-1}ds\\
                 &=&\Gamma(\alpha/2)-\frac{2}{\alpha}x^{\alpha/2}(t+1)^{\alpha/2}\int_0^1
                 e^{-w^{2/\alpha}x(t+1)}dw\/.
                  \end{eqnarray*}
                 Thus we get
                  \begin{eqnarray*}
                     H(x,y) &=& {\sin(\pi\alpha/2)\o\pi}
                \int_1^\infty
                {e^{-(y-x)t}\o(t^2-1)^{\alpha/2}} \, dt\/ -  {2\sin(\pi\alpha/2)x^{\alpha/2}\o\alpha\pi\Gamma(\alpha/2)}
                \int_1^\infty\int_0^1
                 e^{-w^{2/\alpha}x(t+1)}dw
                {e^{-(y-x)t}\o(t-1)^{\alpha/2}} \, dt\/\\
                &=&U_1(x,y)-{2\sin(\pi\alpha/2)x^{\alpha/2}\o\alpha\pi\Gamma(\alpha/2)}
                \int_1^\infty\int_0^1
                 e^{-w^{2/\alpha}x(t+1)}dw
                {e^{-(y-x)t}\o(t-1)^{\alpha/2}}\,dt\/.
                  \end{eqnarray*}
                  Because
                  \begin{equation*}
                  G_{(0,\infty)}(x,y) = U_1(x,y)-
                 E^x[ e^{-\tau_{(0,\infty)}} U_1(X_{\tau_{(0,\infty)}},y)]\/,
                  \end{equation*}
                  thus we finally have
                  \begin{eqnarray*}
                  G_{(0,\infty)}(x,y)&=&
                    {2x^{\alpha/2}\o\alpha\Gamma(\alpha/2)^2\Gamma(1-\alpha/2)}
                  \int_1^\infty\int_0^1e^{-xw^{2/\alpha}(t+1)}\,dw
                  {e^{-(y-x)t}\o(t-1)^{\alpha/2}} \, dt\\
                  &=&
                  {2x^{\alpha/2}\o\alpha\Gamma(\alpha/2)^2\Gamma(1-\alpha/2)}
                  \int_0^1e^{-xw^{2/\alpha}}\int_1^\infty e^{-xtw^{2/\alpha}} {e^{-(y-x)t}\o(t-1)^{\alpha/2}} \,
                  dt\,dw\\
                  &=&
                  {2x^{\alpha/2}e^{x-y}\o\alpha\Gamma(\alpha/2)^2\Gamma(1-\alpha/2)}
                  \int_0^1e^{-2xw^{2/\alpha}}\int_0^\infty e^{-u(xw^{2/\alpha}+y-x)} {du\o u^{\alpha/2}} \,dw\\
                  &=&
                  {2x^{\alpha/2}e^{x-y}\o\alpha\Gamma(\alpha/2)^2}
                  \int_0^1e^{-2xw^{2/\alpha}}(xw^{2/\alpha}+y-x)^{\alpha/2-1}\,dw\\
                  &=&
                  {e^{x-y}\o\Gamma(\alpha/2)^2}
                  \int_0^xe^{-2v}v^{\alpha/2-1}(v+y-x)^{\alpha/2-1}\,dv\\
                  &=&
                  {1\o2^\alpha\Gamma(\alpha/2)^2}
                  \int_0^{4xy}e^{-(u+(y-x)^2)^{1/2}}u^{\alpha/2-1}(u+(y-x)^2)^{-1/2}\,du\\
                  &=&{(y-x)^{{\alpha-1}} \o  2^\alpha \Gamma(\alpha/2)^2}
               \int_0^{{4xy\o(y-x)^2}}e^{(y-x)(t+1)^{1/2}}t^{\alpha/2-1}(t+1)^{-1/2}\,dt\/.
                  \end{eqnarray*}
                 \end{proof}

 %%%%%%%%%%%%%%%%%%%%%%%%%%%%%%%%%%%%%%%%%%%%%%%%%%%%%%%%%%%%%%%%%%%%%%%%%%%%%%%%%%%%%%%%%%%
 %%%%%%%%%%%%%%%%%%%%%%%%%%%%%%%%%%%%%%%%%%%%%%%%%%%%%%%%%%%%%%%%%%%%%%%%%%%%%%%%%%%%%%%%%%%

%%%%%%%%%%%%%%%%%%%%%%%%%%%%%%%%%%%%%%%%%%%%%%%%%%%%%%%%%%%%%%%%%%%%%%%%%%%%%%%%%%%%%%%%%%%
%%%%%%%%%%%%%%%%%%%%%%%%%%%%%%%%%%%%%%%%%%%%%%%%%%%%%%%%%%%%%%%%%%%%%%%%%%%%%%%%%%%%%%%%%%%

 \section{Poisson kernel and Green function of $-(m^{2/\alpha}I-\Delta)^{\alpha/2}$ for $\H$}

 Our proof relies on the computations of the $(d-1)$-dimensional Fourier transform of
 the $d$-dimensional Poisson kernel as
   well as the corresponding $d$-dimensional Green function. Namely we show that these
    Fourier transforms can be expressed in terms of the corresponding
    one-dimensional objects, which we can easily invert.

To avoid confusion, we introduce the following notation, to distinguish one-dimensional
and $d$-dimensional objects.

Notation: $$\text{\bf{x}}\in \R^{d-1},\quad x=(\text{\bf{x}},x_d)\in \R^{d}.$$
  In the proofs below one-dimensional quantities play an important role,
  hence we denote by
 \begin{eqnarray*}&&q^{{m}}_{t}(x)-\mbox{one-dimensional\
  $\alpha$-stable  relativistic density with parameter ${m}$},\\
 &&V^{{m}}_{\lambda}(x)-\mbox{corresponding  $\lambda$-potential},\\
 &&Q^m_{(0,\infty)}(x,u)-\mbox{corresponding  Poisson kernel.}
 \end{eqnarray*}
  We begin with computing of the $(d-1)$-dimensional Fourier transform of the
   $\lambda$-potential $U_{\lambda}((\text{\bf{x}},x_d),(\text{\bf{y}},y_d))=
 U_{\lambda}((\text{\bf{x}}-\text{\bf{y}}),(x_d-y_d))$ with respect to the variable
 $ \text{\bf{y}}$. To avoid notational complexity, we denote it by
 $ \widehat{U_{\lambda}(x,\cdot)}(\text{\bf{z}})$. We employ this kind of notation
 throughout the whole section. Thus, we have

 \begin{lem}
$$ \widehat{U_{\lambda}(x,\cdot)}(\text{\bf{z}})= e^{i(\text{\bf{z}},\text{\bf{x}})}V^{\kappa^\alpha}_{\tilde{\lambda}}(x_d-y_d),
 $$
 where $\kappa= (|\text{\bf{z}}|^2+1)^{1/2}$and $\tilde{\lambda}=\kappa^\alpha+\lambda-1$. Specifying to the case
 $\lambda=1$:
 \begin{eqnarray}
 \widehat{U_1(x,\cdot)}(\text{\bf{z}})
 &=&e^{i(\text{\bf{z}},\text{\bf{x}})}V^{\kappa^\alpha}_{\kappa^\alpha}(x_d-y_d)\nonumber \\
 &=&{2^{1-\alpha \o 2} \o \sqrt{\pi} \Gamma(\alpha/2)}
 e^{i(\text{\bf{z}},\text{\bf{x}})}\,\kappa^{1-\alpha \o 2}
   { K_{1-\alpha \o 2}(\kappa|x_d-y_d|) \o |x_d-y_d|^{1-\alpha \o 2}}
    \label{mpot}.
 \end{eqnarray}
 \end{lem}

 \begin{proof}
 We begin with computation of
 the $(d-1)$-dimensional Fourier transform of the transition density of the  normal distribution:
 $$\widehat{g_u(x,\cdot)}(\text{\bf{z}})=\int_{\R^{d-1}}g_u(x-y)e^{i(\text{\bf{z}},
 \text{\bf{y}})}d\text{\bf{y}}=
 {e^{i(\text{\bf{z}},\text{\bf{x}})} \o {(4\pi u)^{1/2}}}e^{-|\text{\bf{z}}|^2u}
 e^{-\frac {(x_d-y_d)^2}{4u}}.$$
 %%%%%%%%%%%%
 In the next step we use this to find the $(d-1)$-dimensional Fourier
   transform of the transition density $p_t$:
 %%%%%%%%%%%%%%%%%%
 \begin{eqnarray*}
 \widehat{p_t(x,\cdot)}(\text{\bf{z}})&=&
 \int_{\R^{d-1}}p_t(x-y)e^{i(\text{\bf{z}},\text{\bf{y}})}d\text{\bf{y}}
 =e^{t}\int_0^\infty  \hat{g}_u(x)(\text{\bf{z}})e^{-u}\,\theta_{\alpha}(t,u)du\\
 &=&e^{i(\text{\bf{z}},\text{\bf{x}})} e^{t}\int_0^\infty {1 \o {(4\pi u)^{1/2}}}
 e^{-(|\text{\bf{z}}|^2+1)u}  e^{-\frac {(x_d-y_d)^2}{4u}}\theta_{\alpha}(t,u)du\\
 &=&e^{i(\text{\bf{z}},\text{\bf{x}})}
  e^{t(1-(|\text{\bf{z}}|^2+1)^{\alpha/2})}\int_0^\infty {1 \o {(4\pi u)^{1/2}}}
  e^{t(|\text{\bf{z}}|^2+1)^{\alpha/2}}e^{-(|\text{\bf{z}}|^2+1)u}
   e^{-\frac {(x_d-y_d)^2}{4u}}\theta_{\alpha}(t,u)du\\
  &=&e^{i(\text{\bf{z}},\text{\bf{x}})}
  e^{-(\tilde{m}-1)t}\int_0^\infty e^{\tilde{m}t}e^{-\tilde{m}^{2 / \alpha}u}{1 \o
  (4\pi u)^{1/2}} e^{-\frac {(x_d-y_d)^2}{4u}}\theta_{\alpha}(t,u)du
   ,
 \end{eqnarray*}
 where $\tilde{m}= (|\text{\bf{z}}|^2+1)^{\alpha/2}$.
  Note that the the integral expression in the last line is the one-dimensional
  $\alpha$-stable  relativistic density with parameter $\tilde{m}$ which
  we denote by $q^{\tilde{m}}_{t}$. Hence we obtain
 $$ \widehat{p_t(x,\cdot)}(\text{\bf{z}})= e^{i(\text{\bf{z}},\text{\bf{x}})}e^{(1-\tilde{m})t}q^{\tilde{m}}_{t}(x_d-y_d). $$
As a consequence we obtain the
 $(d-1)$-dimensional Fourier  transform of the  $\lambda$-potential :

 \begin{eqnarray*}
 \widehat{U_{\lambda}(x,\cdot)}(\text{\bf{z}})
 &=&e^{i(\text{\bf{z}},\text{\bf{x}})}\int_0^\infty
 e^{(1-\tilde{m}-\lambda)t}q^{\tilde{m}}_{t}(x_d-y_d)dt\\
  &=&e^{i(\text{\bf{z}},\text{\bf{x}})}V^{\tilde{m}}_{(\tilde{m}+\lambda-1)}(x_d-y_d)=
  e^{i(\text{\bf{z}},\text{\bf{x}})}V^{\tilde{m}}_{\tilde{\lambda}}(x_d-y_d),
 \end{eqnarray*}
 where $\tilde{\lambda}=(\tilde{m}+\lambda-1)$ and $V^{\tilde{m}}_{\tilde{\lambda}} $
  is the $\tilde{\lambda} $  potential for the one-dimensional relativistic semigroup
  $q^{\tilde{m}}_{t}$ with parameter $\tilde{m}$.
 Now we take $\lambda=1$ then $\tilde{\lambda}=(\tilde{m}+\lambda-1)=
 (|\text{\bf{z}}|^2+1)^{\alpha/2}=\tilde{m}$. Denote $\kappa=\kappa(z)=
 (|\text{\bf{z}}|^2+1)^{1/2}$. By (\ref{m-potential}) we have
  \begin{eqnarray*}
  V^{\tilde{m}}_{\tilde{\lambda}}(x_d-y_d)&=& V^{\kappa^\alpha}_{\kappa^\alpha}(x_d-y_d)=
  {2^{1-\alpha \o 2} \o \sqrt{\pi} \Gamma(\alpha/2)}
  \kappa^{1-\alpha \o 2}
   { K_{1-\alpha \o 2}(\kappa|x_d-y_d|)\o |x_d-y_d|^{1-\alpha \o 2}}  \,.
  \end{eqnarray*}
\end{proof}

Our next step is to find the $(d-1)$-dimensional Fourier transform  of
 the Poisson kernel and the Green function.

Recall that
 $$P_{\H}(x,u)=E^x[ e^{-\tau_\H},X_{\tau_\H}\in du],\quad x\in\H,\ u\in\H^c\,,$$
 is the Poisson kernel for the operator $-(I - \Delta)^{\alpha/2}$ on the
 set $\H$.
 We also recall the formula \pref{green_def} for the Green function, applied for $m=1$:

 %the transition density of the process killed
 %when leaving the set $\H$:
 %\begin{equation*}
 % p_t^{\H}(x,y)= p_t(x-z)-E^x[\tau_{\H}\le t; p_{t-\tau_{\H}}(X_{\tau_{\H}}-y)]\,.
 %\end{equation*}
 %The Green function is then the density function of the $1$-resolvent operator
 %for the sub-probability semigroup of the killed process. Thus, we obtain
  \begin{equation*}
    G_{\H}(x,y)=
    U_1(x-y)-E^x[ e^{-\tau_\H} U_1(X_{\tau_{\H}},y)],
    \end{equation*}
   for $ x, y \in \H$. The function $G_{\H}(x,y)$ makes sense also for $ x \in \H$ and
    $ y \in \H^c$ and from the general theory it is zero for this range of $x$ and $y$.
    This actually gives a condition for the Poisson kernel:
   %%%%%%%%%%%%%%%%%%%%
   \begin{eqnarray}
   0= G_{\H}(x,y)
   &=&U_1(x-y)-E^x[ e^{-\tau_\H} U_1(X_{\tau_{\H}},y)]\nonumber\\
  &=&U_1(x-y)-\int_{\H^c} P_{\H}(x,u)U_1(u,y)du,\label{rGreen1}
    \end{eqnarray}
for $ x \in \H$ and $ y \in \H^c$.
Now we take the Fourier transform of $G_{\H}(x,y)$. We carry out  computation for
 $x=({\bf0},x_d)$ since the general case can be easily deduced by translation
 invariance
 of the process.    Denote $\kappa=\kappa(z)=
 (|\text{\bf{z}}|^2+1)^{1/2}$ and recall that
 $V^{m}_{m}$ is the $m$-potential for the one-dimensional relativistic semigroup with
 parameter $m$.
 % with the convention that we drop $m$ if $m=1$ .
 %Note that by (\ref{scaling}) we have
 %$V^{m}_{m}(x)= m^{1/\alpha -1}V(m^{1/\alpha}x)$.
   Then,
 \begin{eqnarray*}
  \widehat{G_{\H}(x,\cdot)}(\text{\bf{z}}) &=&
   \widehat{U_1(x,\cdot)}(\text{\bf{z}})  -E^x[ e^{-\tau_\H} \widehat{U_1(X_{\tau_{\H}},\cdot)}(\text{\bf{z}})]\\
  &=& V^{\kappa^\alpha}_{\kappa^\alpha}(x_d-y_d)-E^x[ e^{-\tau_\H}e^{i(\text{\bf{z}},\text{\bf{X}}_{\tau_{\H}})}V^{\kappa^\alpha}_{\kappa^\alpha}(X_{\tau_{\H}}^d-y_d)].
 \end{eqnarray*}
%%%%%%%%%%%%%%%%%%%%%%%%%%%%%%%%%
 %From symmetry considerations   the 1-Poisson kernel  has to have the following form:
 %$$   P_1(x,u)=P_1(x_d,\text{\bf{x}}-\text{\bf{u}},u_d ),  u_d<0<x_d.$$
 %%%%%%%%%%%%%%%%%%%%%%%%%%%%%%%%%%%%%55
 %That is it depends on the the first $d-1$ variables only trough $\text{\bf{x}}-\text{\bf{u}}$. Using this property we %obtain:
 %%%%%%%%%%%%%%%%%%%%%%%%
 %and its Fourier transform (\ref{F_P}) we evaluate:
 %%%%%%%%%%%%%%%%%%%%%
 Next,
 \begin{eqnarray}
 E^x[ e^{-\tau_\H}e^{i(\text{\bf{z}},\text{\bf{X}}_{\tau_{\H}})}
 V^{\kappa^\alpha}_{\kappa^\alpha}(X_{\tau_{\H}}^d-y_d)]
 &=& \int_{u_d<0}e^{i(\text{\bf{z}},\text{\bf{u}})}
 P_{\H}(x,u)V^{\kappa^\alpha}_{\kappa^\alpha}(u_d-y_d)du\nonumber\\
&=& \int^0_{-\infty}\left\{\int_{R^{d-1}}
e^{i(\text{\bf{z}},\text{\bf{u}})}
P_{\H}((\text{\bf{0}},x_d),(\text{\bf{u}},u_d ))d \text{\bf{u}}\right\}
V^{\kappa^\alpha}_{\kappa^\alpha}(u_d-y_d)du_d\nonumber\\
  &=& \int_{-\infty}^0\widehat{P_{\H}(x,\cdot)}(\text{\bf{z}})
  V^{\kappa^\alpha}_{\kappa^\alpha}(u_d-y_d)du_d.\label{compensator}
 \end{eqnarray}
%%%%%%%%%%%%%%%%%
 Hence the condition (\ref{rGreen1}) is equivalent to:

$$ V^{\kappa^\alpha}_{\kappa^\alpha}(x_d-y_d)=
\int_{-\infty}^0\widehat{P_{\H}(x,\cdot)}(\text{\bf{z}})V^{\kappa^\alpha}_{\kappa^\alpha}(u_d-y_d)du_d,
$$
for $y_d <0<x_d$.
 Now this implies by the one-dimensional result ( Theorem \ref{1dPoiss}) that

 \begin{equation}\label{QFourier}\widehat{P_{\H}(x,\cdot)}(\text{\bf{z}})=
Q_{(0,\infty)}^{\kappa^\alpha}(x_d,u_d)=C_{\alpha}^1
 \left({x_d \o -u_d} \right)^{\alpha/2} \,
{e^{-\kappa|x_d-u_d|} \o |x_d-u_d|}.
 \end{equation}
 Note that we used the fact that $\widehat{P_{\H}(x,\cdot)}(\text{\bf{z}})$ is
  real by symmetry of the process.

 Next we proceed to  computation of the Fourier transform for $G_{\H}(x,y)$
 for $ x=({\bf0},x_d), y \in \H$.
 It immediately follows from   (\ref{compensator}), (\ref{QFourier}) and
 Theorem \ref{green_1_form} that
 \begin{eqnarray*}\widehat{G_{\H}(x,\cdot)}(\text{\bf{z}})&=&
  \left(V^{\kappa^\alpha}_{\kappa^\alpha}(x_d-y_d)-
  \int_{-\infty}^0Q_{\H}^{\kappa^\alpha}(x_d,u_d)
  V^{\kappa^\alpha}_{\kappa^\alpha}(u_d-y_d)du_d\right)\\
  &=& G_{(0,\infty)}^{\kappa^\alpha}( x_d, y_d)\\
  &=& {1\o 2^{\alpha} \Gamma(\alpha/2)^2}
  \int_0^{4x_d y_d}
   { s^{{\alpha \o 2}-1} \o (s+(x_d-y_d)^2)^{1/2}}
   e^{-(s+(x_d-y_d)^2)^{1/2}(|\text{\bf{z}}|^2+1)^{1/2} } \,ds\,.
  \end{eqnarray*}

\begin{thm}[Poisson kernel] \label{relPoisson}
%%%%%%%%%%%%%%%%%%%%%%%%
 Let
%%%%%%%%%%%%%%%%%%%%%%%%%%%%%%%%%%
 $$E^x[ e^{-m\tau_\H},X_{\tau_\H}^m\in du]=P_{\H}^m(x,u)\,.$$
%%%%%%%%%%%%%%%%%%%%%%%%%%%%%%%%%%%%%%%%%%%%%%%%%%%%%%%%%%
Then we have
 \begin{equation} \label{dp1}
 P_{\H}^m(x,u)=2C_{\alpha}^1\({m^{1/\alpha} \o 2\pi}\)^{d/2}
 \left({x_d \o -u_d} \right)^{\alpha/2} \,
{ K_{d/2}(m^{1/\alpha}|x-u|) \o |x-u|^{d/2}},\quad
 u_d<0<x_d.
 \end{equation}
\end{thm}
\begin{proof} We prove only the case $m=1$ since the general case follows from the
scaling property. Also it is enough to consider $ x=({\bf0},x_d)$. Applying (\ref{QFourier}) we have
 \begin{eqnarray*}\widehat{P_{\H}(x,\cdot)}(\text{\bf{z}})&=&
 Q_{(0,\infty)}^{\kappa^\alpha}(x_d,u_d)\\
 &=&
C_{\alpha}^1
 \left({x_d \o -u_d} \right)^{\alpha/2} \,
{e^{-\kappa|x_d-u_d|} \o |x_d-u_d|}\\
&=&
C_{\alpha}^1
 \left({x_d \o -u_d} \right)^{\alpha/2} \,
{e^{-(|\text{\bf{z}}|^2+1)^{1/2}|x_d-u_d|} \o |x_d-u_d|}.
\end{eqnarray*}
Taking into account (\ref{relFourier}) and \pref{Cauchyrel} with $(d-1)$ instead of $d$ and $|x_d-y_d|$
instead of $t$ we  complete the proof.
%Then we invert this using (\ref{Cauchyrel}) and (\ref{relFourier}).
\end{proof}

\begin{thm}[Green function for $\H$] \label{relPoisson1}
\begin{equation}
 G_{\H}^m(x,y) = {2^{1-\alpha}m^{d/2\alpha} |x-y|^{\alpha-d/2} \o (2\pi)^{d/2} \Gamma(\alpha/2)^2}
 \int_0^{4 x_d y_d \o |x-y|^2 } {t^{{\alpha \o 2} -1} \o (t+1)^{d/4}}
  K_{d/2}(m^{1/\alpha}|x-y|(t+1)^{1/2})  \, dt\,.
 \end{equation}
 \end{thm}
 \begin{proof}We prove only the case $m=1$ since the general case follows from the scaling property. Also it is enough to consider $ x=({\bf0},x_d)$. Recall that
 $$\widehat{G_{\H}(x,\cdot)}(\text{\bf{z}})= {1 \o 2^{\alpha} \Gamma(\alpha/2)^2}
  \int_0^{4x_d y_d}
   { s^{{\alpha \o 2}-1} \o (s+(x_d-y_d)^2)^{1/2}}
   e^{-(s+(x_d-y_d)^2)^{1/2}(|\text{\bf{z}}|^2+1)^{1/2} } \,ds\,.
 $$
Taking into account (\ref{relFourier}) and \pref{Cauchyrel} with $(d-1)$ instead of $d$ and $(s+(x_d-y_d)^2)^{1/2}$
instead of $t$ we obtain for $d>1$
 \begin{eqnarray*}
 G_{\H}(x,y) &=&  {2^{1-\alpha} \o (2\pi)^{d/2} \Gamma(\alpha/2)^2}
 \int_0^{4 x_d y_d} s^{{\alpha \o 2} -1}
  {K_{d/2}((|\text{\bf{y}}|^2 +(x_d-y_d)^2+s)^{1/2}) \o
 (|\text{\bf{y}}|^2 +(x_d-y_d)^2+s)^{1/2}) ^{d/4}} \, ds \\
  &=&
 {2^{1-\alpha} \o (2\pi)^{d/2} \Gamma(\alpha/2)^2}
 \int_0^{4 x_d y_d} s^{{\alpha \o 2} -1}
 {K_{d/2}((|x-y|^2+s)^{1/2}) \o (|x-y|^2+s)^{d/4}} \, ds \\
 &=& {2^{1-\alpha} |x-y|^{\alpha-d/2} \o (2\pi)^{d/2} \Gamma(\alpha/2)^2}
 \int_0^{4 x_d y_d \o |x-y|^2 } {t^{{\alpha \o 2} -1} \o (t+1)^{d/4}}
  K_{d/2}(|x-y|(t+1)^{1/2})  \, dt\,.
 \end{eqnarray*}
 \end{proof}

 \begin{cor}\label{trivalbound}

 Assume that $|x-y|\le 1$. Then there is $ C = C(\alpha,d)$ such that
$$\left(x_d\,  y_d \wedge1 \right)^{\alpha/2}\le C G_{\H}(x,y).$$
\end{cor}
This estimate is not optimal but sufficient for our purposes in the next section.
%Very precise estimates based on the formulas obtained above are proved in \cite{BRB}.

  \begin{proof} First, observe that $ K_{d/2}(r)/r^{d/2}$ is decreasing
  (see, e.g. \cite{E1}).
  Next, we use one of the form of the Green function obtained in the above proof to arrive at
 \begin{eqnarray*}
 G_{\H}(x,y) &=&
 {2^{1-\alpha} \o (2\pi)^{d/2} \Gamma(\alpha/2)^2}
 \int_0^{4 x_d y_d} s^{{\alpha \o 2} -1}
 {K_{d/2}((|x-y|^2+s)^{1/2}) \o (|x-y|^2+s)^{d/4}} \, ds \\
 &\ge&{2^{1-\alpha} \o (2\pi)^{d/2} \Gamma(\alpha/2)^2}
 \int_0^{1\wedge 4 x_d y_d} s^{{\alpha \o 2} -1}
 {K_{d/2}((|x-y|^2+s)^{1/2}) \o (|x-y|^2+s)^{d/4}} \, ds\\
 &\ge&{2^{1-\alpha} \o (2\pi)^{d/2} \Gamma(\alpha/2)^2}
 \int_0^{1\wedge 4 x_d y_d} s^{{\alpha \o 2} -1}
 {K_{d/2}(2) \o 2^{d/4}} \, ds\\
 &=& C \left(4x_d \, y_d \wedge1 \right)^{\alpha/2}.
 \end{eqnarray*}

 \end{proof}
 Note that using the asymptotic behaviour of $K_{d/2}(r)$ as $r\to 0^+$ we obtain
 the well-known formulas for the Green function and the Poisson kernel for the standard
 symmetric (rotation invariant) $\alpha$-stable process:
 \begin{equation} \label{stableGreen}
 \lim_{m\to 0^+}G_{\H}(x,y) = C(\alpha,d)
 \int_0^{4 x_d y_d \o |x-y|^2 } {t^{{\alpha \o 2} -1} \o (t+1)^{d/2}}
 \, dt, \quad
 x_d,\ y_d>0.
 \end{equation}

 Similarly,

  \begin{equation} \label{stablePoisson}
 \lim_{m\to 0^+}P_{\H}^m(x,u)= C(\alpha,d)
 \left({x_d \o -u_d} \right)^{\alpha/2} \,
{ 1 \o |x-u|^{d}},\quad
 u_d<0<x_d.
 \end{equation}
 As far as we know, it is the first alternative proof of the Poisson kernel formula for
 the standard symmetric (rotation invariant) $d$-dimensional $\alpha$-stable
 process without  application of Kelvin's
 transform.

 As a corollary we compute now the formula for $E^z e^{-\tau_\H}$. Its importance, among
 other things, is due to the fact that it is harmonic for the operator
  $(I-\Delta)^{\alpha/2}$. In probabilistic terms, it is harmonic
 for the Schr\"odinger operator based on the generator of our relativistic process
 with the potential $q=-1$.

 \begin{cor} \label{rtauH}
 For any $0<\alpha<2$ we have
 \begin{equation*}
 E^z e^{-\tau_\H}= { 1 \o {\Gamma(\alpha/2)}} \int_{z_d}^{\infty}t^{\alpha/2-1} e^{-t}\,dt.
  \end{equation*}
 \end{cor}
  \begin{proof}
  The proof consists of computing the mass of the Poisson kernel. It is obvious that we may assume that $d=1$.

Substituting $(-u)=v^{2\o\alpha-2}$ and taking into account the following well-known
    identity
      \begin{eqnarray*}
       {1 \o v^{2\o2-\alpha} + x}
        &=&  \int_0^{\infty} e^{-w(x+v^{2\o2-\alpha})}\,dw\
       \end{eqnarray*}
   we have, after changing  order of
    integration:
    \begin{eqnarray*}
    E^z e^{-\tau_\H}&=& {\sin(\alpha \pi/2) \o \pi}\int_{-\infty}^0 \({z \o -u}\)^{\alpha/2} {e^{-(z-u)} \o z-u}\, du\\
     &=&
    {2\sin(\pi\alpha/2)\o(2-\alpha)\pi}
      z^{\alpha/2}e^{-z} \int_0^\infty \{\int_0^\infty
       e^{-wz}e^{-wv^{2\o2-\alpha}}\,dw\}
       e^{-v^{2\o2-\alpha}} \, dt\,dv
     \\
     &=&
      z^{\alpha/2}e^{-z}\int_0^\infty
       e^{-wz}\int_0^\infty e^{-v^{2\o2-\alpha}(w+1)}\,dv
        \, dw \/.\\
    &=&
   z^{\alpha/2}
    \int_0^\infty
       (w+1)^{\alpha/2-1}e^{-(w+1)z}dw\\
    &=&
    { 2 \o {\Gamma(\alpha/2)}} \int_{\sqrt{z}}^{\infty}t^{\alpha-1} e^{-t^2}\,dt \\
    &=&
    { 1 \o {\Gamma(\alpha/2)}} \int_{z}^{\infty}t^{\alpha/2-1} e^{-t}\,dt
    \end{eqnarray*}
    This clearly ends the proof.

  \end{proof}
  %%%%%%%%%%%%%%%%%%%%%%%%%%%%%%%%%%%%%%%%%%%%%%%%%%%%%%%%%%%%%%%%%%%%%%%%%%%%%%%%%%%%%
  %%%%%%%%%%%%%%%%%%%%%%%%%%%%%%%%%%%%%%%%%%%%%%%%%%%%%%%%%%%%%%%%%%%%%%%%%%%%%%%%%%%%%
 \section{Estimates of Green function of $I-(I-\Delta)^{\alpha/2}$ for $\H$}
In this section
we work under the assumption that $m=1$. In this section we deal with
 the Green function for the set $\H$, corresponding to the operator
 $I-(I-\Delta)^{\alpha/2}$. Equivalently, it is the usual Green function (that is,
 $0$-Green function) for the relativistic process $X_t$. To distinguish it from
 the previously considered Green function for the operator $-(I-\Delta)^{\alpha/2}$
   we denote it  by $G_{\H}^0$.  Our  objective is to establish some estimates
   for $G_{\H}^0$. The estimates will be sharp if $x,y$ are
   close enough. We apply  for this purpose the estimate for  $G_{\H}$ contained in
   Corollary \ref{trivalbound}.

 The following lemma will be useful in the sequel.

\begin{lem}    \label{density} There is $C$ such that:
$$ p_t(x-y) - p_t(x-y^*)\le p_t^{\H}(x,y)
\le  C (t^{-d/2}+ t^{-d/\alpha}) P^x(\tau_\H\ge t/3)P^y(\tau_\H\ge t/3),$$
where $y^*= (y_1, \dots, y_{d-1}, -y_d)$, $x, y \in \H$.
 \end{lem}
\begin{proof} We start with the upper bound.  Since $p_t^{\H}(x,y)$ is a density of a
semigroup and  $p_t^{\H}(x,y)\le \max_{z \in {\R}^d} p_t(z)$ then we have
%%%%%%%%%%%%
$$p_{2t}^{\H}(x,y)=\int_\H p_t^{\H}(x,z)p_t^{\H}(z,y)dz\le
\max_{z \in {\R}^d} p_t(z)\int_\H p_t^{\H}(x,z)dz=
\max_{z \in {\R}^d} p_t(z)P^x(\tau_\H\ge t).$$
%%%%%%%%%%%%%%%%%%
Next we repeat that argument to have
%%%%%%%%%%%%%%%%
\begin{eqnarray*}p_{3t}^{\H}(x,y)&=&\int_\H p_{2t}^{\H}(x,z)p_t^{\H}(z,y)dz\le
\max_{z \in {\R}^d} p_t(z)P^x(\tau_\H\ge t)\int_\H p_t^{\H}(z,y)dz\\
&=&\max_{z \in {\R}^d} p_t(z)P^x(\tau_\H\ge t)P^y(\tau_\H\ge t),
\end{eqnarray*}
which proves the upper bound since
 $\max_{z \in {\R}^d} p_t(z)\le
  C (t^{-d/2}+ t^{-d/\alpha})$ (see Lemma \ref{transden_0}).

To get the lower bound we use the subordination of the process  to the Brownian motion:
$X_t=B_{T_{\alpha}(t)}$.
%%%%%%%%%%%%%
Then $$p_t^{\H}(x,y)=P^x( B_{T_{\alpha}(t)}\in dy,
 B_{T_{\alpha}(s)}\in\H, 0 \le s<t )\ge
P^x( B_{T_{\alpha}(t)}\in dy, B_s\in\H, 0 \le s<T_{\alpha}(t) )$$
%%%%%%%%%%%%%%%%%%%%%
Using the independence of $T_{\alpha}$ and $B$ we obtain
%%%%%%%%%%%%%%%%%%
$$P^x( B_{T_{\alpha}(t)}\in dy, B_s\in\H, 0 \le s<T_{\alpha}(t)|T_{\alpha}(\cdot) )=
 g_{T_{\alpha}(t)}(x-y) - g_{T_{\alpha}(t)}(x-y^*),$$
by the well known formula for the density of the killed Brownian motion  on exiting $\H$.
Integrating and using the fact that $E^0 g_{T_{\alpha}(t)}(z)= p_t(z)$ we obtain the
lower bound.
\end{proof}

\begin{lem}    \label{potential_lower} There is $C$ such that:
$$  G^0_\H(x,y)\ge  C G_{\H}^g(x,y),$$
where  $G_{\H}^g(x,y)$ is the Green function of $\H$ for the Brownian motion.
 \end{lem}
\begin{proof} Let $V(x,y)=\int_0^\infty (p_{t}(x-y)-p_{t}(x-y^*))dt$. From the previous
 lemma it is enough to prove that $V(x,y)\ge  C G_{\H}^g(x,y).$ We have

\begin{eqnarray*}
V(x,y)&=&\int_0^\infty (p_t(x-y)-p_t(x-y^*))dt\\
&=&\int_0^\infty e^{t}\int_0^\infty (g_u(x-y)-g_u(x-y^*))e^{-u}\theta_\alpha(t,u) du dt\\
&=&\int_0^\infty (g_u(x-y)-g_u(x-y^*))e^{-u}\int_0^\infty e^{t}\theta_\alpha(t,u)dt du\\
&=&\int_0^\infty (g_u(x-y)-g_u(x-y^*))G(u) du,
\end{eqnarray*}
where  $G(u)=e^{-u}\int_0^\infty e^{t}\theta_\alpha(t,u)dt$ is the potential
 of the subordinator $T_\alpha(t)$. It was proved in \cite{RSV} that $G(u)$ is
  a completely monotone (hence decreasing) function and
$\inf_{u> 0}G(u)=\lim_{u\to \infty}G(u)=C$. %To show that it is enough to find that  Laplace transfor $\lim_{u\to \infty}G(u)=C>0$.
 We find the constant $C=\lim_{u\to \infty}G(u)$ by taking into account
 the asymptotics of the
  Laplace transform of $G(u)$ at the origin:
%%%%%%%%%%%%%%
\begin{eqnarray*}
\int_0^\infty e^{-\lambda u}G(u) du
= \int_0^\infty e^{t}\int_0^\infty e^{-u(1+\lambda)}\theta_\alpha(t,u)du dt
= \int_0^\infty e^{t}e^{-(1+\lambda)^{\alpha/2} t}dt=
 \frac 1 {(1+\lambda)^{\alpha/2} -1}\sim \frac 2{\lambda \alpha}
\end{eqnarray*}
Applying the Tauberian Theorem we obtain that $C= 2/ \alpha$.

 Thus, since $g_u(x-y)-g_u(x-y^*)\ge 0$ we finally obtain
\begin{eqnarray*}
V(x,y)
&=&\int_0^\infty (g_u(x-y)-g_u(x-y^*))G(u) du\\
&\ge&C \int_0^\infty (g_u(x-y)-g_u(x-y^*)) du = C G_{\H}^g(x,y).\\
\end{eqnarray*}
\end{proof}

The next result provides a general upper bound for the Green function.
 \begin{thm}    \label{upperG_0}
 There is a constant $C$ such that
 %%%%%%%%%%%%%%%%%%%%%%%%%%%%
 \begin{equation} \label{u_G_0}
 G_\H^0(x,y)\le  C [\psi(x_d)\psi(y_d)+G_{\H}(x,y)], \quad x,y\in \H,
  \end{equation}
 where  the function $\psi $ has the following form depending on the dimension $d$:

 for $d=1$: $\psi(v)= v^{\alpha/2},\quad 0<v<1$ and $\psi(v)= v^{1/2} ,\quad v\ge 1\,;$

 for $d=2$: $\psi(v)= v^{\alpha/2},\quad 0<v<1$ and $\psi(v)= \ln^{1/2} \(1+ v\) ,\quad v\ge 1\,;$

 for $d\ge 3$: $\psi(v)= (v\wedge 1)^{\alpha/2},\quad v>0\,.$
 \end{thm}
\begin{proof}

The proof will rely on estimates of $P^x(\tau_\H\ge t)$ and the application of
Lemma \ref{density}.

 We proceed to estimate the Green function.
First we split the  integration
$$\int_0^\infty p_t^{\H}(x,y)dt =  \int_0^{2} p_t^{\H}(x,y)dt +
\int_{2}^{\infty} p_t^{\H}(x,y)dt.$$
%%%%%%%
To estimate the first integral we use
%%%%%%%%%%
%$$p_\H(t,x,x)\le C \frac 1{t^{1/\alpha}}, \quad t\le 1$$
%%%%%%%%%%%%
$$\int_0^{2} p_t^{\H}(x,y)dt\le e^2\int_0^{2}e^{-t} p_t^{\H}(x,y)dt
 \le e^2\, G_{\H}(x,y).$$
%%%%%%%%%%%%
Due to Lemma  \ref{density}
$$\int_{2}^{\infty} p_t^{\H}(x,y)dt\le \int_{2}^{\infty}
 P^x(\tau_\H\ge t)P^y(\tau_\H\ge t)\, {dt \o t^{d/2}}=R(x,y). $$

For $d\ge3$ we have
%%%%%%%%%%%%%%
$$R(x,y)\le P^x(\tau_\H\ge 2)P^y(\tau_\H\ge 2)\int_{2}^{\infty}
{dt \o t^{d/2}}\le C(x_d\wedge 1)^{\alpha/2} (y_d\wedge 1)^{\alpha/2}\,, $$
due to (\ref{uboundtail}), which completes the proof in this case.

To deal with $d=1, 2$ observe that by the Schwarz inequality
$$R^2(x,y)\le R(x,x)R(y,y),$$
so we estimate $R(x,x)$. First consider the case $x_d\le \sqrt{2}$.
Then using \pref{Grzywny1} and \pref{uboundtail} we have
%%%%%%%%%%%%%%%%
$$R(x,x)=\int_{2}^{\infty} \left(P^x(\tau_\H\ge t)\right)^2 \,{dt \o t^{d/2}}
\le C x_d^{\alpha}
\int_{2}^{\infty}  \frac {(\ln t)^2}{t} {dt \o t^{d/2}}.$$
%%%%%%%%%%%%%%%
If  $x_d>\sqrt{2}$, using (\ref{Grzywny1}) we estimate \begin{eqnarray*}
R(x,x)=\int_{2}^{\infty}\left(P^x(\tau_\H\ge t)\right)^2\, {dt \o t^{d/2}}
&=&
 \int_{2}^{x_d^2} \left(P^x(\tau_\H\ge t)\right)^2\, {dt \o t^{d/2}}
 +\int_{x_d^2}^\infty \left(P^x(\tau_\H\ge t)\right)^2 {dt \o t^{d/2}}\\
&\le& C\int_1^{x_d^2}\, {dt \o t^{d/2}} +
C\int_{x_d^2}^\infty \left(\frac{x_d + \ln t}{t^{1/2}}\right)^2\, {dt \o t^{d/2}}.%\le C y,
\end{eqnarray*}
Hence for  $d=1$ we have  $R(x,x)\le C x$, while for  $d=2$ we arrive at $R(x,x)\le C \ln x_2$.
Taking into account all cases  we get
%%%%%%%%%%%%%%
$$\int_{2}^\infty p_t^{\H}(x,y)dt \le C \psi(x)\psi(y).$$
The proof
of the theorem is complete.

\end{proof}

%%%%%%%%%%%%%%%%%%%%%%%%%%%%%%%%%%%%%
%%%%%%%%%%%%%%%%%%%%%%%%%%%%%%%%%%%%
\begin{thm} \label{bounds}

  For $d=  1$ and  $|x-y|<1, $ we have that
 % \begin{equation}\label{gbound1}
 %G_{\H}^{r,1}(x,y)+ C_1 (x\wedge y) \le G_{\H}^{r,0}(x,y) \le C [G_{\H}^{r,1}(x,y)+x\wedge y]\,.
 %\end{equation}

 \begin{equation}\label{gbound1}
G_{\H}^0(x,y) \approx   G_{\H}(x,y)+  x\wedge y\,.
 \end{equation}

 For $d= 2 $ and  $|x-y|<1$ we have that
  \begin{equation}\label{gbound2}
  G_{\H}^0(x,y) \approx G_{\H}(x,y)+ \ln(1\vee(x_2\wedge y_2))\,.
 \end{equation}

 For $d\ge 3$ we have that for $|x-y|<1$,
   \begin{equation}\label{gbound3}  G_{\H}^0(x,y) \approx  G_{\H}(x,y).\end{equation}

 \end{thm}

% All the estimates from above follow directly from the Theorem \ref{upperG_0}
% and appriopriate estimates of $G_{\H}^1(x,y)$.

 \begin{proof}

Step {\bf 1}. We first show that if $|x-y|<1$ then

  \begin{equation} \label{Green0_d}
  G_{\H}(x,y) \approx G_{\H}^0(x,y)\, \quad {\text{if}} \quad x_d \wedge y_d \leq 2
  \quad  {\text{or if}} \quad d \geq 3.
  \end{equation}
  Recall that by  Corollary \ref{trivalbound} we have
  \begin{equation} \label{trivial}
  (x_d\,y_d\wedge 1)^{\alpha/2}\le C G_{\H}(x,y).
\end{equation}
  We first consider the case
    $x_d\wedge y_d\le 2$. The  condition $|x-y|<1$ implies that $x_d\vee y_d\le 3$. Hence
    $$x_d\, y_d \le 9(1\wedge x_d\, y_d).$$
   Then due to Theorem \ref{upperG_0} and (\ref{trivial}) we obtain
   %%%%%%%%%%%%%%%%%%%%%%%%
   \begin{equation} \label{Green0_d0}
   G_{\H}(x,y) \le G_{\H}^0(x,y)
   \le C_1(G_{\H}(x,y)+ (x_d\, y_d)^{\alpha/2}) \le C_2G_{\H}(x,y)\,.
    \end{equation}
  %%%%%%%%%%%%%%%%%%%%%%%%%%%%%%%%%%%%%%%%

 We now consider the case when $d\ge 3$.
 %The  condition $|x-y|<1$ implies that $|x_d-y_d|<1$
 %and this yields
 It is elementary that
 $$(x_d \wedge 1)(y_d \wedge 1)\le x_d\, y_d \wedge 1.$$
 Then the inequality \pref{Green0_d0}
 can be rewritten
 in the following way regardless on the assumption on $x_d\wedge y_d$,
 %%%%%%%%%%%%%%%%%%%%%%%
 \begin{eqnarray*} %\label{Green0_d1}
 G_{\H}(x,y) &\le& G_{\H}^0(x,y)\le C(G_{\H}(x,y)+
  ((x_d \wedge 1)(y_d \wedge 1))^{\alpha/2})\\
  &\le& C(G_{\H}(x,y)+
  (x_d\,y_d \wedge 1)^{\alpha/2})\\
  &\le&
 C_1G_{\H}(x,y),
 \end{eqnarray*}
 where in the last line we applied again  (\ref{trivial}).
 This completes the proof of Step {\bf{1}} and
 proves the theorem for $d\ge 3$.

 Step {\bf{2}}.
 In this step we complete the proof of the case $d=1$.

 The lower bound follows from the estimate proved in Lemma \ref{potential_lower} :
 $G_{(0,\infty)}^0(x,y)\ge C G_{(0,\infty)}^g(x,y)$, where $G_{(0,\infty)}^g$
 is the Green function of
 $(0,\infty)$ for the Brownian motion.

 It is well known that $G_{(0,\infty)}^g(x,y)= x\wedge y$ so
\begin{equation*}
  x\wedge y +
  G_{(0,\infty)}(x,y) \leq C G_{(0,\infty)}^0(x,y)\,,
\end{equation*}
and the lower bound is established. In  the case $x\wedge y \le 2$ we also get
the upper bound since
 $ G_{(0,\infty)}^0(x,y) \approx  G_{(0,\infty)}(x,y))$ by Step {\bf 1}.
 %application of Theorem \ref{upperG_0} yields

If $x\wedge y \ge 2$ we obtain by  Theorem \ref{upperG_0},
 %%%%%%%%%%%%%%%%%%%%%%
 \begin{equation} \label{Green0_1}
  G_{(0,\infty)}^0(x,y)
 \le C_2(G_{(0,\infty)}(x,y)+
  (x y)^{1/2}).
 % \le C_3(G_{(0,\infty)}^1(x,y)+ x\wedge y).
 \end{equation}
 %%%%%%%%%%%%%%%%%%%
  Since $|x-y|<1$ and $x\wedge y \geq 2$  we have
  $(xy)^{1/2} \leq x \vee y\le   (x \wedge y)+1\leq  2 (x \wedge y) $,
   which completes the proof of
 \pref{gbound1} for the case $x\wedge y \geq 2$.

 Step {\bf{3}}.
 Now  we deal with the case $d=2$.
 We claim that  $G_{\H}^g(x,y) \geq C \ln (1\vee (x_2\wedge y_2))$. It is enough to
 show it for
 $x_2\wedge y_2\ge 1$. It is well known that
 $G_{\H}^g(x,y)= {1 \o 2\pi}\ln\frac {|x^*-y|}{|x-y|}$.
 If $|x-y|\le 1$ then %also $y_2-x_2<1$ so
 %%%%%%%%%%%%%%
  \begin{eqnarray*}
  G_{\H}^g(x,y)&=& {1 \o 2\pi}\ln\frac {|x^*-y|}{|x-y|}
  \ge {1 \o 2\pi}\ln{|x^*-y|}\\
  &\ge& {1 \o 2\pi}\ln {(x_2+y_2)} \geq {1 \o 2\pi} \ln (1\vee (x_2\wedge y_2)).
   \end{eqnarray*}
  As a consequence of the above inequality and Lemma \ref{potential_lower} we obtain
  \begin{equation} \label{basic_estim_d2}
        G_{\H}(x,y) \le  G_{\H}(x,y) + \ln(1\vee (x_2\wedge y_2))
     \le C( G_{\H}(x,y) +  G_{\H}^g(x,y)) \le  G_{\H}^0(x,y)
    \end{equation}

  If now $ x_2\wedge y_2\le 2$ then by Step {\bf{1}}
  $ G_{\H}(x,y) \approx   G_{\H}^0(x,y)$ and we obtain the conclusion.

  We proceed with the case $ x_2\wedge y_2\ge 2$.
   Note that $x_2\vee y_2\le 1+ x_2\wedge y_2$, so
 $$ \ln^{1/2} (1+ y_2) \ln^{1/2}(1+ x_2)\le 2 \ln (x_2\wedge y_2)=
  2 \ln (1\vee (x_2\wedge y_2)) $$
  Hence  by Theorem \ref{upperG_0},
  $$ G_{\H}^0(x,y) \le C (G_{\H}(x,y) + \ln^{1/2} (1+ y_2) \ln^{1/2}(1+ x_2))
  \le
  C (G_{\H}(x,y)+ \ln (1\vee (x_2\wedge y_2))),$$
  which, together with the estimate \pref{basic_estim_d2} completes the proof
 of Step {\bf{3}} and of the theorem.

 \end{proof}
%%%%%%%%%%%%%%%%%%%%%%%%%%%%%%%%%%%%%%
%%%%%%%%%%%%%%%%%%%%%%%%%%%%%%%

 {\bf Remarks.}

 {\bf 1.} The above theorem shows that, in particular, for fixed $y \in \H$ the function
 $G_{\H}^0(x,y) \approx x_d^{\alpha/2} $, for $x$ near the boundary. With an extra effort
 we are able to strengthen this result and actually show that
 for fixed $y\in \H$ we have
 \begin{equation*}
\lim_{x \to \xi\in \partial\H} G_{\H}^0(x,y)/ x_d^{\alpha/2}=C(y,\xi)\,.
 \end{equation*}

 {\bf 2.} We are able to show that for $d=1$ the estimate for $ G_{\H}^0$ stated in Theorem
 \ref{bounds} is optimal for the whole range of $x, y \in \H$, that is, we
 have
 \begin{equation*}
  G_{\H}^0(x,y) \approx  G_{\H}(x,y) + x\wedge y\,, \quad x\,, y \in \H\,.
 \end{equation*}

\end{document}